\documentclass[11pt]{article}

\usepackage{amsfonts}
\usepackage{shuffle,yfonts,amsmath,amssymb,amsthm}
\usepackage{hyperref}
\hypersetup{pdfstartview={XYZ null null 1.00}}
\usepackage{pdfpages}
\usepackage{cite}
\newcommand\nc{\newcommand}
 \nc{\binn}{{\binom{2n}{n}}}

\newcommand{\cl}{{\rm Cl}}

 \allowdisplaybreaks

\catcode`!=11
\let\!int\int \def\int{\displaystyle\!int}
\let\!lim\lim \def\lim{\displaystyle\!lim}
\let\!sum\sum \def\sum{\displaystyle\!sum}
\let\!sup\sup \def\sup{\displaystyle\!sup}
\let\!inf\inf \def\inf{\displaystyle\!inf}
\let\!cap\cap \def\cap{\displaystyle\!cap}
\let\!max\max \def\max{\displaystyle\!max}
\let\!min\min \def\min{\displaystyle\!min}
\let\!frac\frac \def\frac{\displaystyle\!frac}
\catcode`!=12

\let\oldsection\section
\renewcommand\section{\setcounter{equation}{0}\oldsection}

\allowdisplaybreaks

 \nc{\gam}{{\gamma}}
 \nc{\gG}{{\Gamma}}
 \nc{\vep}{{\varepsilon}}
 \nc{\gs}{{\sigma}}
 \nc{\gth}{{\theta}}
 \nc{\gS}{{\Sigma}}
 \nc{\gf}{{\phi}}
 \nc{\gk}{{\kappa}}
 \nc{\gz}{{\zeta}}
 \nc{\tgz}{{\tilde{\zeta}}}
 \nc{\gO}{{\Omega}}
 \nc{\sif}{{\mathcal S}}
 \nc{\gt}{{\tau}}

\def\ze{\zeta}

\theoremstyle{plain}
\newtheorem{thm}{Theorem}[section]

\theoremstyle{definition}

\begin{document}
\title{\bf A Note on Sun's Conjectures on Ap\'ery-like Sums Involving Lucas Sequences and Harmonic Numbers}
\author{
{Ce Xu${}^{a,}$\thanks{Email: cexu2020@ahnu.edu.cn, first author, ORCID 0000-0002-0059-7420.}\ \ and Jianqiang Zhao${}^{b,}$\thanks{Email: zhaoj@ihes.fr, corresponding author, ORCID 0000-0003-1407-4230.}}\\[1mm]
\small a. School of Mathematics and Statistics, Anhui Normal University, Wuhu 241002, P.R. China\\
\small b. Department of Mathematics, The Bishop's School, La Jolla, CA 92037, USA}

\date{}
\maketitle

\noindent{\bf Abstract.} In this paper, we will prove Zhi-Wei Sun's four conjectural identities
on Ap\'{e}ry-like sums involving Lucas sequences and harmonic numbers by using a few results of Davydychev--Kalmykov \cite{DavydychevDe2004}.

\medskip



\section{Introduction}
The classical Lucas numbers $\{L_n\}_{n\geq 0}$ are defined by
\begin{align*}
L_0=2,\ L_1=1\quad \text{and}\quad L_{n+1}=L_n+L_{n-1} \ \forall n\ge 1.
\end{align*}
Obviously, for any $n\ge 0$,
$$L_n=\Bigl(\frac{1+\sqrt{5}}{2}\Bigr)^n+\Bigl(\frac{1-\sqrt{5}}{2}\Bigr)^n.$$

Let $A$ and $B$ be integers. Similarly to Sun \cite{Sun2009}, we define the more general Lucas sequence $v_n=v_n(A,B)$ as follows:
\begin{align*}
v_0(A,B):=2,\quad v_1(A,B):=A\quad \text{and}\quad v_{n+1}(A,B)=Av_{n}(A,B)-Bv_{n-1}(A,B) \ \forall n\ge 1.
\end{align*}
The characteristic equation $x^2-Ax+B=0$ of the sequence $\{v_n\}_{n\geq 0}$ has two roots
\[\alpha=\frac{A+\sqrt{A^2-4B}}{2}\quad \text{and}\quad \beta=\frac{A-\sqrt{A^2-4B}}{2}.\]
It is well known that for any integer $n\ge 0$ we have
\begin{align*}
v_n(A,B)=\alpha^n+\beta^n.
\end{align*}
For examples, setting $\phi=\frac{\sqrt{5}+1}{2}$ we have
\begin{align*}
&v_n(3,1)=L_{2n}=\phi^{2n}+\phi^{-2n},\\
&v_n(4,1)=(2-\sqrt{3})^n+(2+\sqrt{3})^n,\\
&v_n(4,2)=(2-\sqrt{2})^n+(2+\sqrt{2})^n,\\
&v_n(5,5)=(\sqrt{5}\phi)^n+(\sqrt{5}\phi^{-1})^n.
\end{align*}

Define $H_0:=0$ and the classical harmonic numbers
\[
H_n:=\sum_{k=1}^n \frac{1}{k}\quad\forall n\ge 1.
\]
In this short note, we will prove the following results. 
\begin{thm}\label{thm-sun's-conj} We have
\begin{align} \label{equ:10.62}
&\sum_{n=1}^\infty \frac{L_{2n}}{n^2\binom{2n}{n}}\left(H_{2n}-H_{n-1}\right)=\frac{41\ze(3)+4\pi^2\log(\phi)}{25},\\
&\sum_{n=1}^\infty \frac{v_n(5,5)}{n^2\binom{2n}{n}}\left(H_{2n}-H_{n-1}\right)=\frac{124\ze(3)+\pi^2\log\left(5^5\phi^6\right)}{50},\label{equ:10.63}\\
&\sum_{n=1}^\infty \frac{v_n(4,1)}{n^2\binom{2n}{n}}\left(H_{2n}-H_{n-1}\right)=\frac{23\ze(3)+2\pi^2\log(2+\sqrt{3})}{12},\label{equ:2022a}\\
&\sum_{n=1}^\infty \frac{v_n(4,2)}{n^2\binom{2n}{n}}\left(H_{2n}-H_{n-1}\right)=\frac{259\ze(3)+2\pi^2(4\log(2)+8\log(1+\sqrt{2}))}{128}.\label{equ:2022b}
\end{align}
\end{thm}

Motivated by his study of the congruences in \cite{Sun2009}, Sun \cite{Sun2015,Sun2021,Sun2022} proposes
identities \eqref{equ:10.62} and \eqref{equ:10.63} as \cite[Conjectures 10.62 and 10.63]{Sun2021}
(see also \cite[Eqs. (3.11) and (3.12)]{Sun2015}). He further conjectures \eqref{equ:2022a}
and \eqref{equ:2022b} in  \cite{Sun2022}.

\section{Proof of Theorem \ref{thm-sun's-conj}}
Recall that for any positive integer $m$ the Clausen function $\cl_m(\theta)$ is defined by
\begin{align*}
&{\cl}_{2m-1}(\theta):=\sum_{n=1}^\infty \frac{\cos(n\theta)}{n^{2m-1}}\quad \text{and}\quad {\cl}_{2m}(\theta):=\sum_{n=1}^\infty \frac{\sin(n\theta)}{n^{2m}} \qquad \forall \theta\in [0,\pi].
\end{align*}
Many special values of the  Clausen function are closely related to the Riemann zeta values. For example,
\begin{equation}\label{equ:Riemmann}
\cl_{2m-1}(0)=\zeta(2m-1),\quad \cl_{2m-1}(\pi)=(2^{2-2m}-1)\zeta(2m-1),\quad \cl_{2m-1}\Big(\frac\pi2\Big)=\frac{\cl_{2m-1}(\pi)}{2^{2m-1}}.
\end{equation}
The following relations among special values of the  Clausen function are well known (see, e.g., \cite[Eq. (6.51)]{Lewin1981}):
for any positive integer $r$ we have
\begin{align}
&\cl_3\big(\tfrac{\pi}{r}\big)+\cl_3\big(\tfrac{3\pi}{r}\big)+\cdots+\cl_3\big(\tfrac{(2r-1)\pi}{r}\big)=-\tfrac3{4r^2}\ze(3),\label{equ-cl-fun}\\
&\cl_3\big(\tfrac{2\pi}{r}\big)+\cl_3\big(\tfrac{4\pi}{r}\big)+\cdots+\cl_3\big(\tfrac{(2r-2)\pi}{r}\big)=-\big(1-\tfrac1{r^2}\big)\ze(3).\label{equ-cl-fun-1}
\end{align}
Setting $r=4,5,6$ in \eqref{equ-cl-fun} and $r=5$ in \eqref{equ-cl-fun-1} one obtains
\begin{align}
&\cl_{3}\big(\tfrac{\pi}{4}\big)+\cl_{3}\big(\tfrac{3\pi}{4}\big)=-\tfrac3{128}\ze(3), \qquad
\cl_{3}\big(\tfrac{\pi}{5}\big)+\cl_{3}\big(\tfrac{3\pi}{5}\big)=\tfrac9{25}\ze(3), \label{equ-cl-fun-exa-1}\\
&\cl_{3}\big(\tfrac{\pi}{6}\big)+\cl_{3}\big(\tfrac{5\pi}{6}\big)=\tfrac1{12}\ze(3), \qquad
\cl_{3}\big(\tfrac{2\pi}{5}\big)+\cl_{3}\big(\tfrac{4\pi}{5}\big)=-\tfrac{12}{25}\ze(3),\label{equ-cl-fun-exa-2}
\end{align}
since $\cl_3(\theta)=\cl_3(2\pi\pm \theta)$, $\cl_3(\pi)=-\tfrac3{4}\ze(3)$ and $\cl_3\big(\tfrac{\pi}{2}\big)=-\tfrac3{32}\ze(3)$ by \eqref{equ:Riemmann}.

On the other hand, Davydychev and Kalmykov have shown that
(see \cite[(2.36), (2.37) and (2.67)]{DavydychevDe2004})
\begin{align*}
&\sum_{n=1}^\infty \frac{1}{n^3\binom{2n}{n}}u^n=2\cl_3(\theta)+2\theta \cl_2(\theta)-2\ze(3)+\theta^2\log\Big(2\sin\frac{\theta}{2}\Big), \\
&\sum_{n=1}^\infty \frac{H_{n-1}}{n^2\binom{2n}{n}}u^n=4\cl_3(\pi-\theta)-2\theta\cl_2(\pi-\theta)+3\ze(3), \\
&\sum_{n=1}^\infty \frac{H_{2n-1}}{n^2\binom{2n}{n}}u^n=-2\cl_3(\theta)+4\cl_3(\pi-\theta)-2\theta\cl_2(\pi-\theta)-\theta\cl_2(\theta)+5\ze(3),
\end{align*}
where $u=4\sin^2\frac{\theta}{2}$ for any $\theta\in [0,\pi]$.
Combining the above three equations we get
\begin{align}\label{equ-comb}
\sum_{n=1}^\infty \frac{u^n}{n^2\binom{2n}{n}}\left(H_{2n}-H_{n-1}\right)&=\sum_{n=1}^\infty \frac{u^n}{n^2\binom{2n}{n}}\left(H_{2n-1}-H_{n-1}+\frac1{2n}\right) \nonumber\\
&=\ze(3)-\cl_3(\theta)+\frac{\theta^2}{2}\log\Big(2\sin\frac{\theta}{2}\Big).
\end{align}
Taking special values of $u$ we obtain
\begin{align*}
u=\Big(\frac{\sqrt{5}\pm 1}2\Big)^2=4\sin^2\frac{\theta}{2}\quad \Longrightarrow \quad& \sin\frac{\theta}{2}=\frac{\sqrt{5}\pm 1}{4}\quad \text{and}\quad \theta=\frac{2\pi}{5}\pm\frac{\pi}{5},\\
u=\frac{5\pm \sqrt{5}}{2}=4\sin^2\frac{\theta}{2}\quad \Longrightarrow \quad& \sin\frac{\theta}{2}=\sqrt{\frac{5\pm\sqrt{5}}{8}}\quad \text{and}\quad \theta=\frac{3\pi}{5}\pm\frac{\pi}{5},\\
u=2\pm\sqrt{3}=4\sin^2\frac{\theta}{2}\quad \Longrightarrow \quad&
\sin\frac{\theta}{2}=\frac{\sqrt{2\pm\sqrt{3}}}{2}\quad \text{and}\quad \theta=\frac{\pi}{2}\pm\frac{\pi}{3},\\
u=2\pm\sqrt{2}=4\sin^2\frac{\theta}{2}\quad \Longrightarrow \quad& \sin\frac{\theta}{2}=\frac{\sqrt{2\pm\sqrt{2}}}{2}\quad \text{and}\quad \theta=\frac{\pi}{2}\pm\frac{\pi}{4}.
\end{align*}
Therefore, \eqref{equ-comb} and \eqref{equ-cl-fun-exa-1}--\eqref{equ-cl-fun-exa-2} yield
\begin{align*}
\sum_{n=1}^\infty \frac{L_{2n}}{n^2\binom{2n}{n}}\left(H_{2n}-H_{n-1}\right)
&=2\ze(3)-\cl_{3}\big(\tfrac{3\pi}{5}\big)-\cl_{3}\big(\tfrac{\pi}{5}\big)+\frac{4\pi^2}{25}\log(\phi) \\
&=\frac{41\ze(3)+4\pi^2\log(\phi)}{25},\\
\sum_{n=1}^\infty \frac{v_n(5,5)}{n^2\binom{2n}{n}}\left(H_{2n}-H_{n-1}\right)
&=2\ze(3)-\cl_{3}\big(\tfrac{2\pi}{5}\big)-\cl_3\big(\tfrac{4\pi}{5}\big)+\frac{\pi^2\log\left(5^5\phi^6\right)}{50}\\
&=\frac{124\ze(3)+\pi^2\log\left(5^5\phi^6\right)}{50},\\
\sum_{n=1}^\infty \frac{v_n(4,1)}{n^2\binom{2n}{n}}\left(H_{2n}-H_{n-1}\right)
&=2\ze(3)-\cl_{3}\big(\tfrac{\pi}{6}\big)-\cl_{3}\big(\tfrac{5\pi}{6}\big)+\frac{\pi^2\log(2+\sqrt{3})}{6}\\
&=\frac{23\ze(3)+2\pi^2\log(2+\sqrt{3})}{12},\\
\sum_{n=1}^\infty \frac{v_n(4,2)}{n^2\binom{2n}{n}}\left(H_{2n}-H_{n-1}\right)
&=2\ze(3)-\cl_{3}\big(\tfrac{\pi}{4}\big)-\cl_{3}\big(\tfrac{3\pi}{4}\big)+\frac{\pi^2(4\log(2)+8\log(1+\sqrt{2}))}{64}\\
&=\frac{259\ze(3)+2\pi^2(4\log(2)+8\log(1+\sqrt{2}))}{128}.
\end{align*}

This concludes the proof of Theorem \ref{thm-sun's-conj} and this note.

\medskip

\noindent{\bf Acknowledgement.}  Ce Xu expresses his deep gratitude to Prof. Zhi-Wei Sun for valuable discussions and comments. Ce Xu is supported by the National Natural Science Foundation of China (Grant No. 12101008), the Natural Science Foundation of Anhui Province (Grant No. 2108085QA01) and the University Natural Science Research Project of Anhui Province (Grant No. KJ2020A0057). Jianqiang Zhao is supported by the Jacobs Prize from The Bishop's School.

\end{document}